\documentclass[twocolumn,showpacs,preprintnumbers,amsmath,amssymb,floatfix,superscriptaddress]{revtex4}
\usepackage{dcolumn}% Align table columns on decimal point
\usepackage{bm}% bold math
\usepackage{color}
\usepackage{amsmath,amsthm,amssymb}
\usepackage[pdftex]{graphicx}
\newcommand{\R}{\mathbb{R}}

\newtheorem{thm}{Theorem}[section]

\newtheorem{prop}{Proposition}
\newtheorem{define}[thm]{Definition}
\newtheorem{remark}[thm]{Remark}
\newcommand{\e}{\epsilon}

\newcommand{\1}{\mathbf{1}}

\newcommand{\av}[1]{\left|{#1}\right|}

\newcommand{\ip}[2]{\left\langle{{#1}},{{#2}}\right\rangle}

\newcommand{\be}{\begin{enumerate}}
\newcommand{\bi}{\begin{itemize}}
\newcommand{\ee}{\end{enumerate}}
\newcommand{\ei}{\end{itemize}}
\newcommand{\ii}{\item}

\newcommand{\E}{\mathbb{E}}

\renewcommand{\span}{{\mathrm{span}}}
\newcommand{\asin}{\sin^{-1}}
\newcommand{\acos}{\cos^{-1}}
\newcommand{\atan}{\tan^{-1}}

\begin{document}

\title{Phase-locked Patterns of the Kuramoto Model on 3-Regular Graphs} \author{Lee
  DeVille} \affiliation{Department of Mathematics, University of
  Illinois} \author{Bard Ermentrout} \affiliation{Department of
  Mathematics, University of Pittsburgh}

\begin{abstract}
  We consider the existence of non-synchronized fixed points to the
  Kuramoto model defined on sparse networks, specifically, networks
  where each vertex has degree exactly three.  We show that ``most''
  such networks support multiple attracting phase-locked solutions
  that are not synchronized, and study the depth and width of the
  basins of attraction of these phase-locked solutions.  We also show
  that it is common in ``large enough'' graphs to find phase-locked
  solutions where one or more of the links has angle difference
  greater than $\pi/2$.
\end{abstract}


\maketitle
\section{Introduction}

\subsection{Background}

Rhythmic phenomena are ubiquitous in biological and physical systems
and are typically modeled as a limit cycle or a network of many
coupled oscillators. Networks of coupled oscillators are, in general,
difficult to analyze, but, in some cases, it is possible to reduce
such networks into similar networks of much simpler models. For
example, if coupling is ``weak'' (in a sense that can be made precise,
see e.g.~\cite{EK84}) and the individual oscillators are nearly
identical, models of the form:
\[
x_i' = f(x_i) + \epsilon \sum_{ij} G_{ij}(x_i,x_j), \quad i,j=1,\ldots,N
\]
where $\epsilon$ is small and $u'=f(u)$ admits an asymptotically
stable limit cycle, can be reduced the the study of a system on the
$N-$torus of the form:
\begin{equation}
\label{eq:general}
\theta_i' = \epsilon \sum_j H_{ij}(\theta_j-\theta_i).
\end{equation}
Here $H_{ij}(\phi)$ are $T-$periodic functions of $\phi$ where $T$ is
the period of the uncoupled limit cycles and
$x_i(t)=u(t+\theta_i)+O(\epsilon).$ When $H_{ij}(\phi) = \omega_i +
K_{ij} \sin(\phi)$, we have the so-called Kuramoto model. The
particular case when $K_{ij}=1/N$ and $N$ is large has been the
subject of much theory~\cite{Strogrev} and serves as a core
model for synchronization phenomena. When $K_{ij}$ is zero except when
$|i-j|=1$, then this represents a network of oscillators in a chain
which has been used as a model for swimming of the lamprey
\cite{cohen-holmes} and other locomotory pattern generators.  More
general functions $H_{ij}(\phi)$ were considered in \cite{Ermentrout.92} where
general conditions for the existence of stable periodic solutions were
given.

In this paper, we want to make a distinction between different types
of synchronization. In the Kuramoto model, synchrony refers to a state
in which there is a collective rhythm or organization that emerges
from the infinite sized system. For a finitely sized system, as we
study here, we are interested in fully {\em phase-locked solutions}
which correspond to periodic solutions to (\ref{eq:general}) and wish
to distinguish them from what we will call {\em synchronous
  solutions}, where $\theta_i(t)=\theta_j(t)$ for all $i,j.$ As we
will formally define below, a phase pattern is a stable phase-locked
state that is not synchronized (in the sense above). For example in a
four oscillator model for quadruped
locomotion~\cite{Golubitsky.Stewart.06} would call the walk (where
each the oscillators take the four phases, $0,1/4,1/2,3/4$) and the
trot (pairs of oscillators are half a cycle out of phase) patterns,
but not the pronk (where all oscillators have the same phase, that is,
they are synchronized in our sense). A more concrete example can be
found in~\cite{paullet1994stable} where we showed that a square array
of nearest neighbor oscillators of the form:
\[
\theta_i' = 1+\sum_{j\in N_i} \sin(\theta_j-\theta_i)
\]   
could admit both a stable synchronized solution ($\theta_i=t$) or
various types of stable rotating waves.  The latter are a pattern and
the former not.  In \cite{udeigwe} we showed that there was a stable
rotating wave on a network of oscillators on the dodecahedral graph
and generalizations of it as well as stable
synchronization. Nonsynchronized activity in the form of phase-waves
is a common feature of spatio-temporal dynamics in the nervous system
\cite{kleinfeld}. This patterned activity can arise from either
heterogeneities in the medium or from topological defects
\cite{paullet1994stable, allen2013} where all the oscillators are
identical.  It is the latter types of patterns that we are interested
in. We finally note that the same network can often admit multiple
stable solutions; such is the case in \cite{paullet1994stable,
  udeigwe, allen2013}.

\subsection{Motivation}

The main object of study in this paper is the Kuramoto system on
sparse graphs, and in particular we concentrate on cubic graphs in
this paper. (Recall that a {\em cubic} graph is a graph which is
regular of degree $3$, i.e. each node has exactly three neighbors in
the graph.)  By Kuramoto system, we mean here that the interaction
between oscillators is the sine-function: in
equation~(\ref{eq:general}) we choose $H_{ij}(\phi)=1+c_{ij}
\sin(\phi)$ where $C=c_{ij}$ is the adjacency matrix for an undirected
graph. In absence of coupling, all oscillators have the same dynamics.
We will use the term ``patterns'' to refer to non-synchronized stable
phase-locked states. The two cases where patterns are well understood
are the cases of ring graphs and of the complete graph.  The ring
graph can support a large number of
patterns~\cite{Wiley.Strogatz.Girvan.06} and the complete graph
supports no patterns.  Moreover, the patterns on the ring are very
easy to characterize.  If $N$ is the number of oscillators in the
ring, then for any $k < \lfloor N /4\rfloor$, there is a stable fixed
point given by the $k$-twist solution $\theta_m = 2\pi k m/N$. The
condition for stability says that the angle difference between any two
connected nodes must be between $-\pi/2$ and $\pi/2.$

Since rings are regular of degree $2$, cubic graphs are the sparsest
regular graphs not yet understood.  We might expect that as we add
nodes, patterns are easier to find as they were in rings.  But of
course the situation is more complex: in the case of the ring, we know
the size of the (only) cycle in the graph, but a cubic graph will have
many cycles and we have to be able to find a set of angles that match
around all of these cycles (note that the angle differences around any
cycle will need to sum to a multiple of $2\pi$).

We are also particularly interested in studying patterns with long
links, by which we mean an edge on which there is an angle difference
larger than $\pi/2$.  As we discuss below, it is straightforward to
show that a pattern that has only short links will be stable, but the
converse is not true.  Of course, it is easy to obtain a long-link
pattern by making the weight on one or more edges particularly weak
and allowing it to stretch.  But it is not {\em a priori} clear that
it is possible to do so when all of the edge weights have the same
magnitude.  It was shown in~\cite{NON-D} that there are no long-link
patterns on the ring when the edge weights are the same, and to the
best of our knowledge, it has never been shown that a Kuramoto model
where the oscillators have identical frequencies and the edge weights
are equal can exhibit a long-link pattern, but in this paper we
exhibit many examples of these.  (The existence of long-link patterns
has been shown on the ring where the edge weights are equal but the
intrinsic frequencies are chosen to be different
in~\cite{Bronski.DeVille.Ferguson.15}.)

\subsection{Definitions}

\begin{define}
  Let $G = (V,E)$ be a graph and define the {\bf Kuramoto flow on $G$}
  as the ODE on $\R^{|V|}$:
  \begin{equation}\label{eq:K}
    \frac{d}{dt} \theta_v = \sum_{w \in N(v)} \sin(\theta_w-\theta_v). 
  \end{equation}
  (Here $N(v)$ is the set of neighbors of the vertex $v$ in the
  graph.)  A {\bf fixed point} or {\bf phase-locked solution} is any
  point $\theta$ which gives a fixed point for this ODE.  We assume
  throughout that the graph $G$ is connected.
\end{define}

\begin{remark}
\bi

\ii We are assume here that all of the oscillators have the
same intrinsic frequency.  It is then true that the strength of
coupling is irrelevant, so we set it to unity.

\ii Since the right-hand side of the ODE is invariant under the
transformation $\theta\mapsto \theta + c\1$, fixed points are only
unique up to translation.  We will typically make the convention of
choosing a distinguished vertex $v$ and enforcing that $\theta_v = 0$.

\ii It is clear from inspection that $\theta_v\equiv 0$ is a fixed
point of~\eqref{eq:K} for any graph $G$.  We will always call this the
{\bf synchronized} or {\bf sync} solution.  Also notice that since the
oscillators are assumed to have the same intrinsic frequencies, this
system will not have precessing solutions with a moving center of
gravity.  Summing~\eqref{eq:K} over $v$ and using the symmetry of the
graph we see that the sum of the angles is constant.

\ei
\end{remark}

\begin{define}
  We say that a fixed point $\theta$ of~\eqref{eq:K} is {\bf stable}
  if the Jacobian at $\theta$ has one zero eigenvalue and the
  remainder negative, and if ${\bf 1}$ is the eigenvector in the
  nullspace.  The reason for this is that, as is easy to show,
  $(d/dt)\ip{\theta}{{\bf 1}} = 0$ for any solution, and thus we can
  pick the sum of $\theta_i$ a priori and see that it is fixed.
  Moreover, the motion in the ${\bf 1}$ direction gives a zero
  eigenvalue, so if it is unique then the point is dynamically stable
  in ${\bf 1}^\perp$.

  Every stable fixed point has an open set which is its basin of
  attraction.  We will refer to the {\bf width} of such a basin as its
  Lebesgue measure on the set ${\bf 1}^\perp$ inside $[0,2\pi]^N$. 
\end{define}

\begin{prop}
  The sync solution is always stable.
\end{prop}

\begin{proof}
  We see that the Jacobian of~\eqref{eq:K} at the sync solution is the
  $|V|\times|V|$ matrix whose entries are given by
\begin{equation*}
  J_{vw} = \begin{cases} 1,&v\neq w, (v,w)\in E,\\ 0,& v\neq w, (v,w)\not\in E,\\ -\deg(v), & v=w.\end{cases}
\end{equation*}
As such it is the negative of the standard graph Laplacian.  As seen
in~\cite{Chung.book}, this implies that $J$ is negative semidefinite,
and the number of zero eigenvalues is the number of connected
components of the graph.  In particular, if $G$ is connected, then the
Jacobian at zero has exactly one zero eigenvalue, and the remainder
are negative.  Moreover, $\ker(J) = \span({\bf 1})$.
\end{proof}

\begin{define}
  Any stable fixed point of~\eqref{eq:K} that is not the sync solution
  is termed a {\bf stable pattern}, or typically just a {\bf pattern}.
\end{define}

\begin{define}
  For any graph $G$, there is a potential energy function
  \begin{equation}\label{eq:defofPhi}
    \Phi_G(\theta) := \sum_{(v,w)\in E(G)} (1-\cos(\theta_v-\theta_w)).
  \end{equation}
  It is not hard to show that~\eqref{eq:K} can be written as $\theta'
  = -\nabla \Phi_G(\theta)$.  The {\bf energy} of a solution is just
  the number $\Phi_G(\theta)$.  Notice that the sync solution always
  has zero energy, and any pattern will have a positive energy.
\end{define}

\begin{define}
  Let $\theta$ be a pattern for some graph $G$.  For each edge
  $(v,w)\in E(G)$, we say the link is {\bf long} if
  $\av{\theta_v-\theta_w} > \pi/2$, {\bf short} if
  $\av{\theta_v-\theta_w} < \pi/2$, and {\bf critical} if
  $\av{\theta_v-\theta_w} = \pi/2$.
\end{define}

\begin{remark}
  The reason for the distinction between short and long is motivated
  by considerations of stability.  Let $\theta$ be a fixed point
  of~\eqref{eq:K} for a given graph $G$.  Then the Jacobian of the
  flow at this point is given by
\begin{equation*}
  J_{vw}= \begin{cases} \cos(\theta_v-\theta_w),&v\neq w, (v,w)\in E,\\ 0,& v\neq w, (v,w)\not\in E,\\ -\sum_{u\in V} J_{vu}, & v=w.\end{cases}
\end{equation*}
Note that each row of this matrix is zero sum, so that $J{\bf 1} =
{\bf 0}$.  Moreover, if all of the off-diagonal terms are positive,
i.e. all of the links are short, then the standard theory of weighted
graph Laplacians~\cite{Chung.book} tells us that the fixed point is
stable.  The case where some of the links are short and some are
critical is addressed in~\cite{Ermentrout.92}.

When there are long links, i.e. when some of the off-diagonal terms
are negative, then there is no assurance that the point is stable.
However, the existence of long links is not enough to exclude
stability: in~\cite{BD-SIAP, BDK-SIMAX}, the study of weighted
Laplacians with signed entries was undertaken and general conditions
were found to determine when such Laplacians would be stable.  In
particular, if there is a configuration with long links such that
every vertex is connected by a path of short links, then it is
possible that this gives a stable point.

\end{remark}

Finally, we note that every thing mentioned above applies to general graphs,
but the subject of the current study is of a restricted class of
graphs, namely:

\begin{define}
  We say a graph is {\bf regular} if every vertex has the same degree,
  and it is {\bf $d$-regular} if that degree is $d$.  A {\bf cubic}
  graph is a 3-regular graph.  It is not hard to see that a
  $d$-regular graph has $\frac12 d|V|$ edges, and so cubic graphs have
  $\frac32 |V|$ edges.  As such, every cubic graph must have an even
  number of vertices.
\end{define}

\subsection{Questions addressed in this paper}\label{sec:questions}

In this paper, we will restrict to only considering {cubic} graphs.
All of the cubic graphs with $N$ vertices have been enumerated up to
$N=40$, and datasets of these graphs are available online,
e.g.~\url{http://staffhome.ecm.uwa.edu.au/~00013890/remote/cubics/}.

Since these data files exist, we can do a comprehensive study of all
cubic graphs of a certain order, and their ability to support
patterns.

The main abstract questions that we address are as follows:

\begin{enumerate}

\item How many cubic graphs support patterns?  A similar question: if
  we choose a cubic graph at random, what are the odds that it
  supports a pattern?

\item Are there graphs that support multiple patterns?  If so, are
  these related?

\item Are there any graphs that support patterns with long links?  If
  yes, is there a graph that supports multiple patterns, some of which
  have long links and some of which do not?  How common are long-link
  patterns?

\item What can we say about the width of the basin of attraction of a
  pattern? Is there a connection between the width of this basin of
  attraction and the energy of a pattern?  Is there a connection
  between this width and the linear stability of a pattern?  Is there
  a connection between width and depth of a basin of attraction?

\end{enumerate}

\subsection{Methods}

In this paper, we are finding attracting fixed points for a system of
ODEs by a Monte Carlo method: we choose random initial conditions
distributed uniformly in $[0,2\pi]^{|V|}$ and flow~\eqref{eq:K} until
the vector field approaches a fixed point (which we detect by
observing the $\ell^2$ norm of the vector field being less than
$10^{-5}$). We then double-check that this is a sink (and not a ``slow
saddle'') by numerically computing the eigenvalues of the Jacobian and
verifying that it is negative semidefinite.  In the simulations, there
were only a tiny number of cases where the numerical Jacobian was
needed because the system was trapped near a slow saddle. For a given
stable pattern, we define the {\em spectral gap} of the pattern as the
distance of the real part of the eigenvalue least negative real part
from zero.  Fixed points with large spectral gaps are in a sense
``more stable'' as perturbations decay more rapidly to the pattern.

\newcommand{\nsamp}{K_{\mathsf{samp}}}

This method both detects attracting fixed points and gives a
statistical estimate of the width of the basin of attraction of a
solution.  For example, if we write $S$ as the basin of a particular
fixed point, then we can determine for any $z\in [0,2\pi]^{|V|}$ if
$z\in S$ or $z\not\in S$.  We then choose $\nsamp$ samples uniformly
in the torus.  The proportion of samples in $S$ will converge to the
normalized measure of $S$ as $\nsamp\to\infty$, and in fact this
random proportion is normally distributed with variance $O(1/\nsamp)$,
so that we could measure this set with whatever confidence we needed.
Also note that this method will never give a false positive for a
pattern, but patterns with small basins of attraction might be missed.
For specific graphs of interest and that do not have too many nodes,
it is possible to find {\em all} fixed points and assess their
stability by using methods from computational algebraic
geometry~\cite{mehta15}. In the future, we hope to use their methods
to study {\em interesting} graphs whose exact fixed points cannot be
readily found.

Finally, we use the enumerated datasets for the cubic graphs at the
link mentioned above.  For small enough $|V|$, we are thus able to
consider every single cubic graph with $|V|$ vertices.  In the last
section of this paper, we discuss several specific graphs that have
some interesting properties; e.g., the smallest cubic graph with a
long link pattern.

\section{Global Statistics}\label{sec:global}

%\subsection{}

First note that we can show, just by exhaustive search, that there are
graphs that support multiple patterns: in fact, we have found one that
supports twelve distinct patterns ($N=30$) and one with ten ($N=20$).

Also, we are careful not to be multiply counting patterns that are
symmetries of each other; for example, any automorphism of $G$ will
move a stable pattern to another, and this one could look quite
different numerically.  However, one can see that any such
automorphism will preserve the energy of a pattern as defined
in~\eqref{eq:defofPhi} --- so we require that two patterns have
different energies before we consider them distinct.

\onecolumngrid

\begin{table*}[ht]
\begin{center}
\begin{tabular}{|c|c|c|c|c|c|c|c|c|c|}\hline
$N\setminus  k$ & total cubics & no patterns & 1 pattern & 2 patterns & 3 patterns &4 patterns&5&6&$> 6$ \\\hline
$10$ & 19 &16& 3  & & &&&&\\\hline
12 & 85 & 61& 22 & 2 &&&&&\\\hline
14 & 509 & 338 & 140 & 27 & 4&&&&\\\hline
16 & 4060 & 2038& 1445 & 457 & 103 & 7&&&\\\hline
18 & 41301 & 17658& 13714 &7048  &2382  &447 &49&2& \\\hline
20 & 9910* & 2533&3565 &2315  &1088  &333 &66&9& 1\\\hline
30 & 1500* &8& 53 & 130 & 245 & 265 & 257& 203& 339\\\hline
\end{tabular}
\caption{Number of graphs supporting a particular number of
  patterns. The last two rows in this table represent a random sample
  of the possible graphs since there are 510,489 graphs for $N=20$ and
  845,480,228,069 graphs for $N=30.$ An exhaustive search would be
  expensive. }
\end{center}
\label{tab:patterns}
\end{table*}
\twocolumngrid

\subsection{Fractions of patterns}

For each $N$, there is $c(N)<\infty$, the total number of cubic graphs with $N$ vertices.  
Let us define $f(N,k)$ as the proportion of these that support $k$ distinct patterns.

%\begin{equation}\label{eq:defoff}
%  f(N,k) = \frac{\#{\mbox{cubic graphs with $N$ vertices that support $k$ distinct patterns}}}{\#{\mbox{cubic graphs wi%th $N$ vertices}}}.
%\end{equation}

\begin{figure}[ht]
\begin{centering}
  \includegraphics[width=0.35\textwidth]{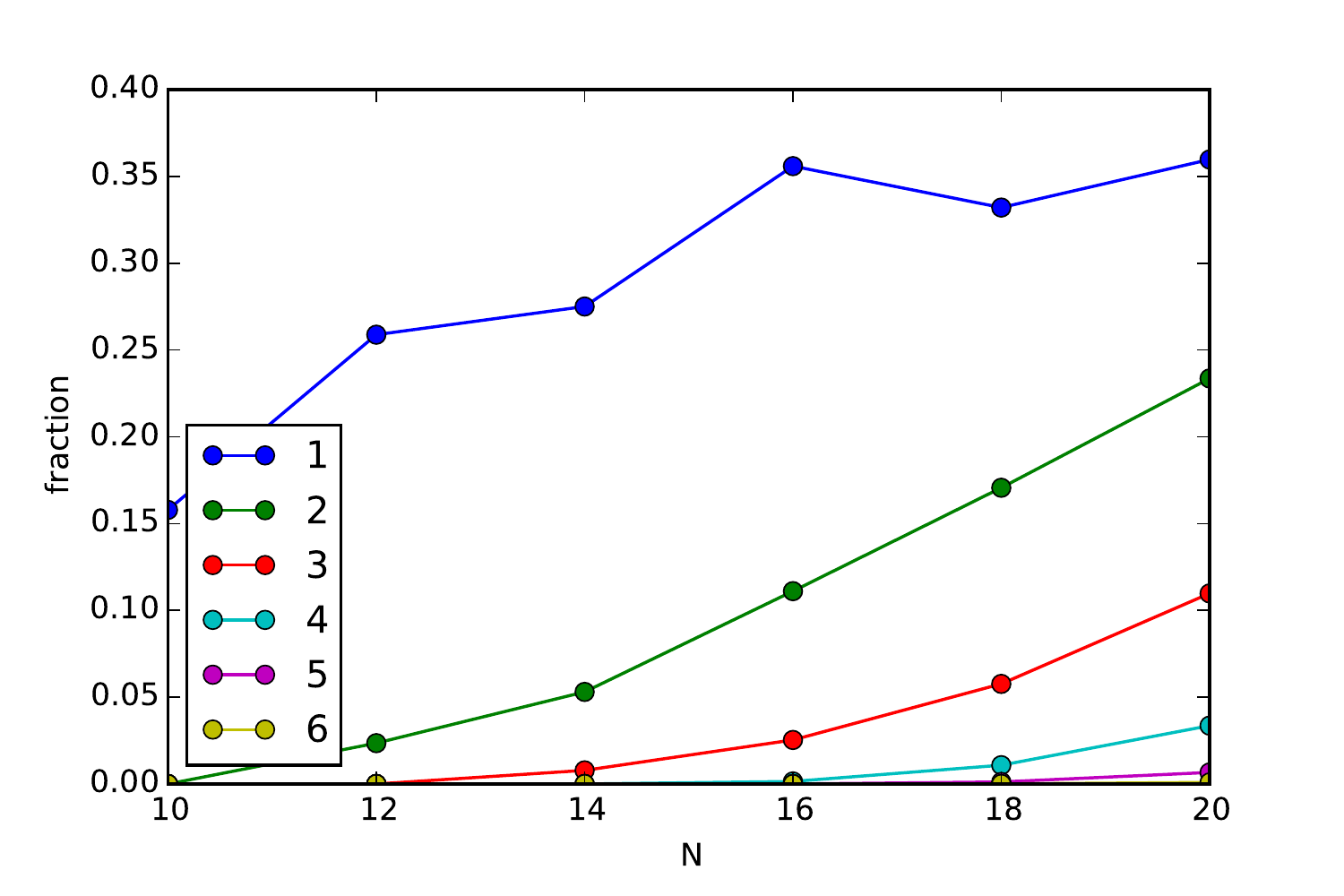}
  \includegraphics[width=0.35\textwidth]{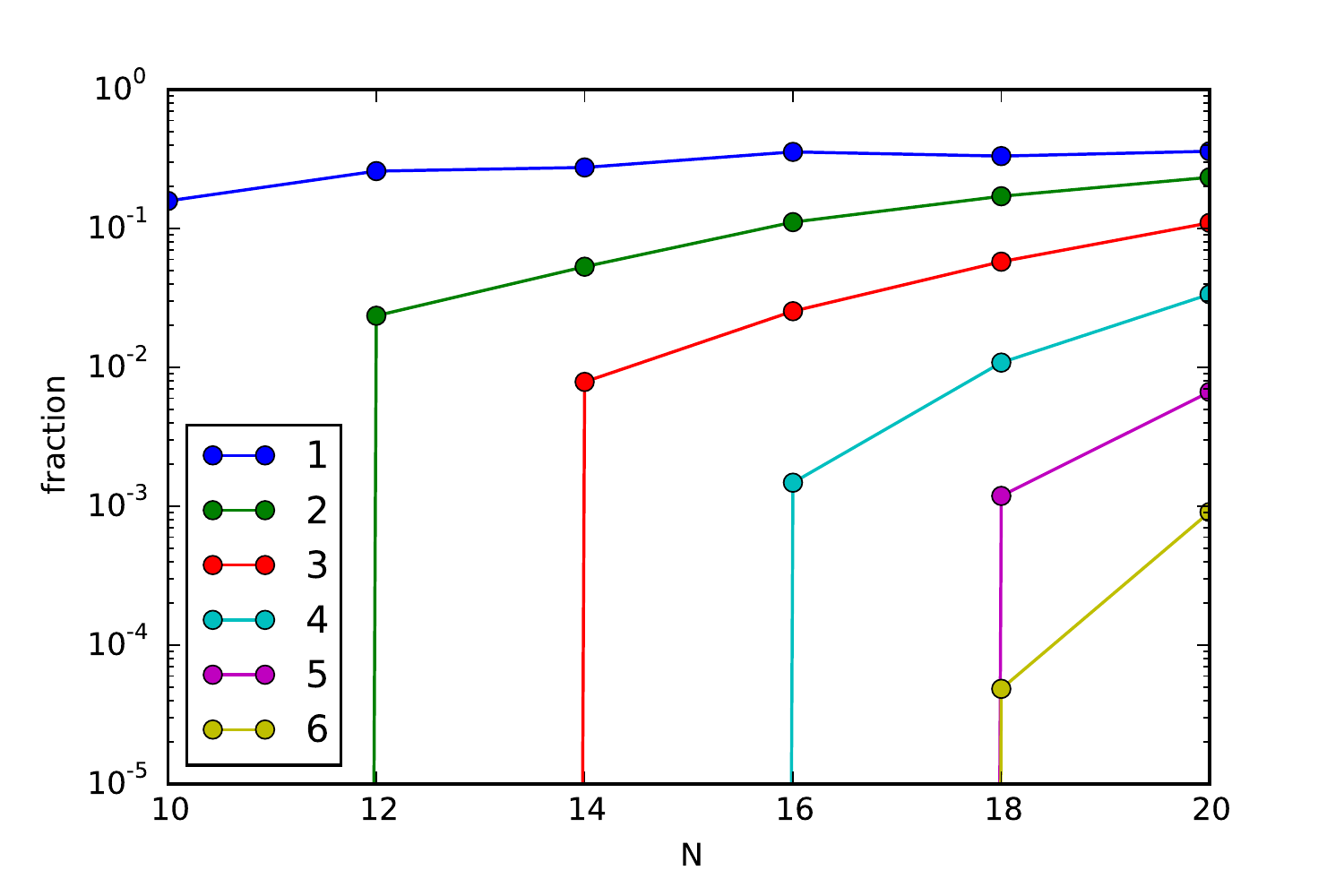}
  \caption{ We plot $f(N,k)$, where different curves correspond to
    different $k$ and we plot versus $N$.  The same data is plotted in
    both frames, one on a linear--linear scale, the other on a
    log--linear scale.}
  \label{fig:pvn}
\end{centering}
\end{figure}

\begin{figure}[ht]
\begin{centering}
  \includegraphics[width=0.35\textwidth]{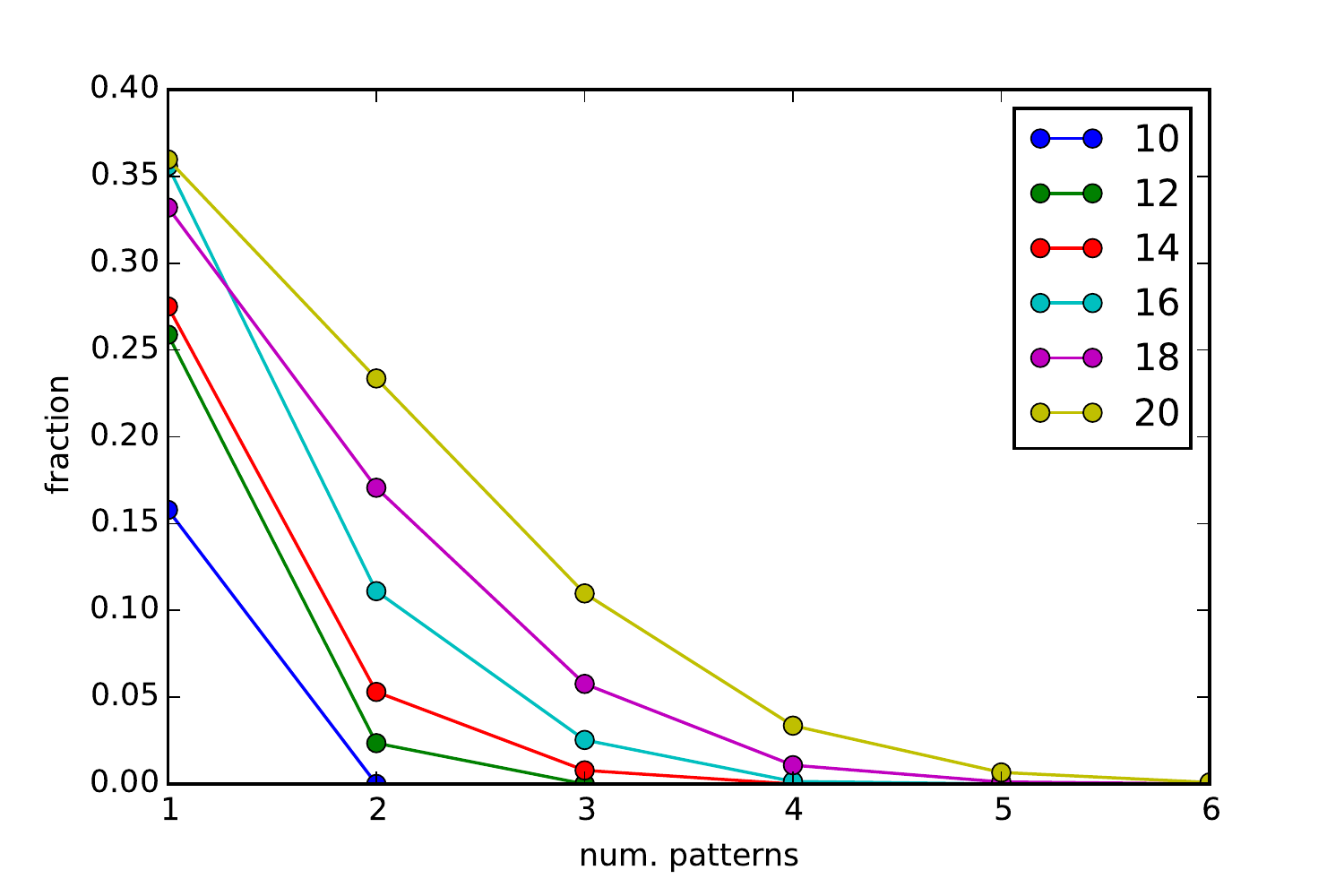}  
  \includegraphics[width=0.35\textwidth]{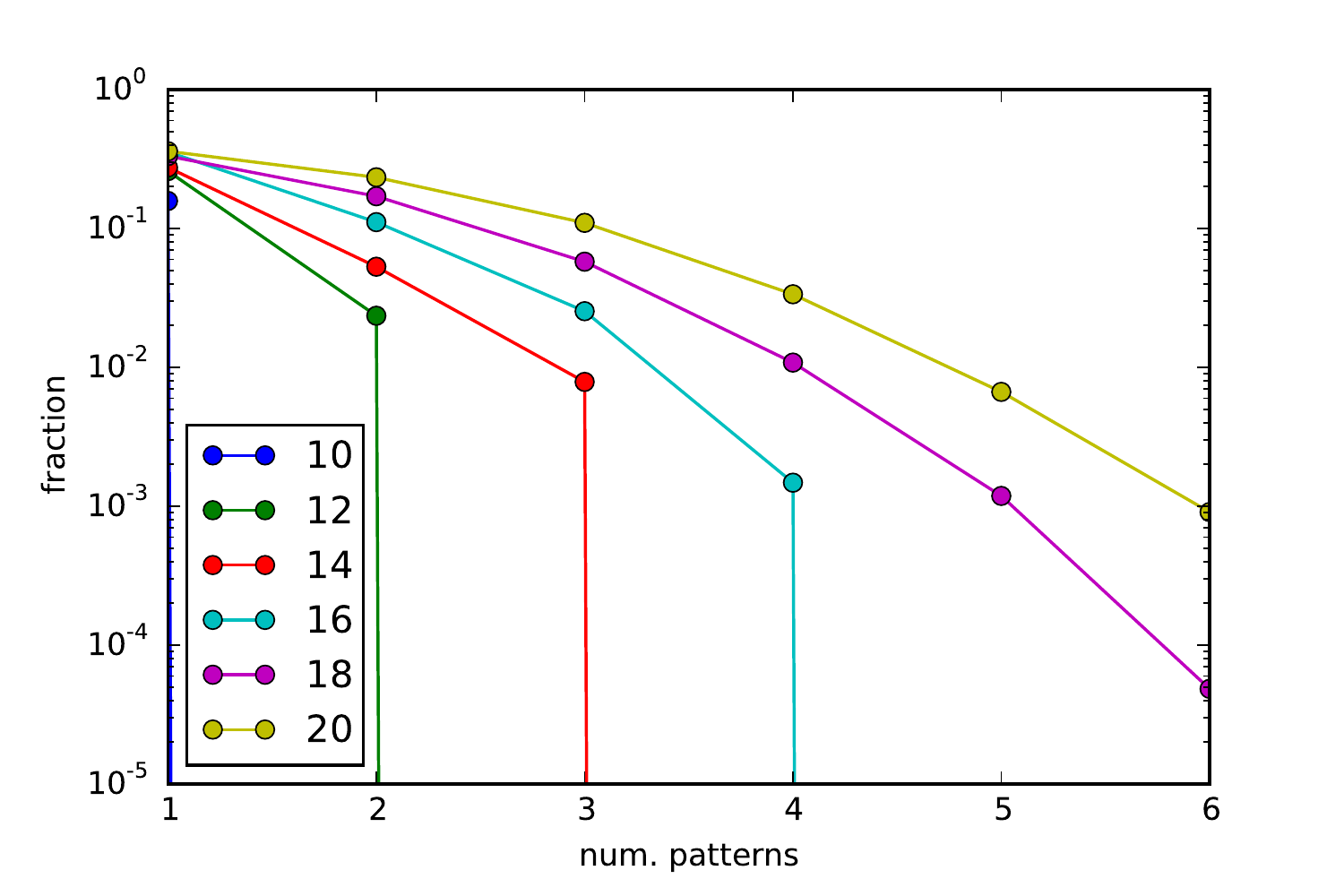}
  \caption{ We plot $f(N,k)$, where different curves correspond to
    different $N$ and we plot versus $k$.  The same data is plotted in
    both frames, one on a linear--linear scale, the other on a
    log--linear scale.}
  \label{fig:pvp}
\end{centering}
\end{figure}

\begin{figure}[ht]
\begin{centering}
  \includegraphics[width=0.5\textwidth]{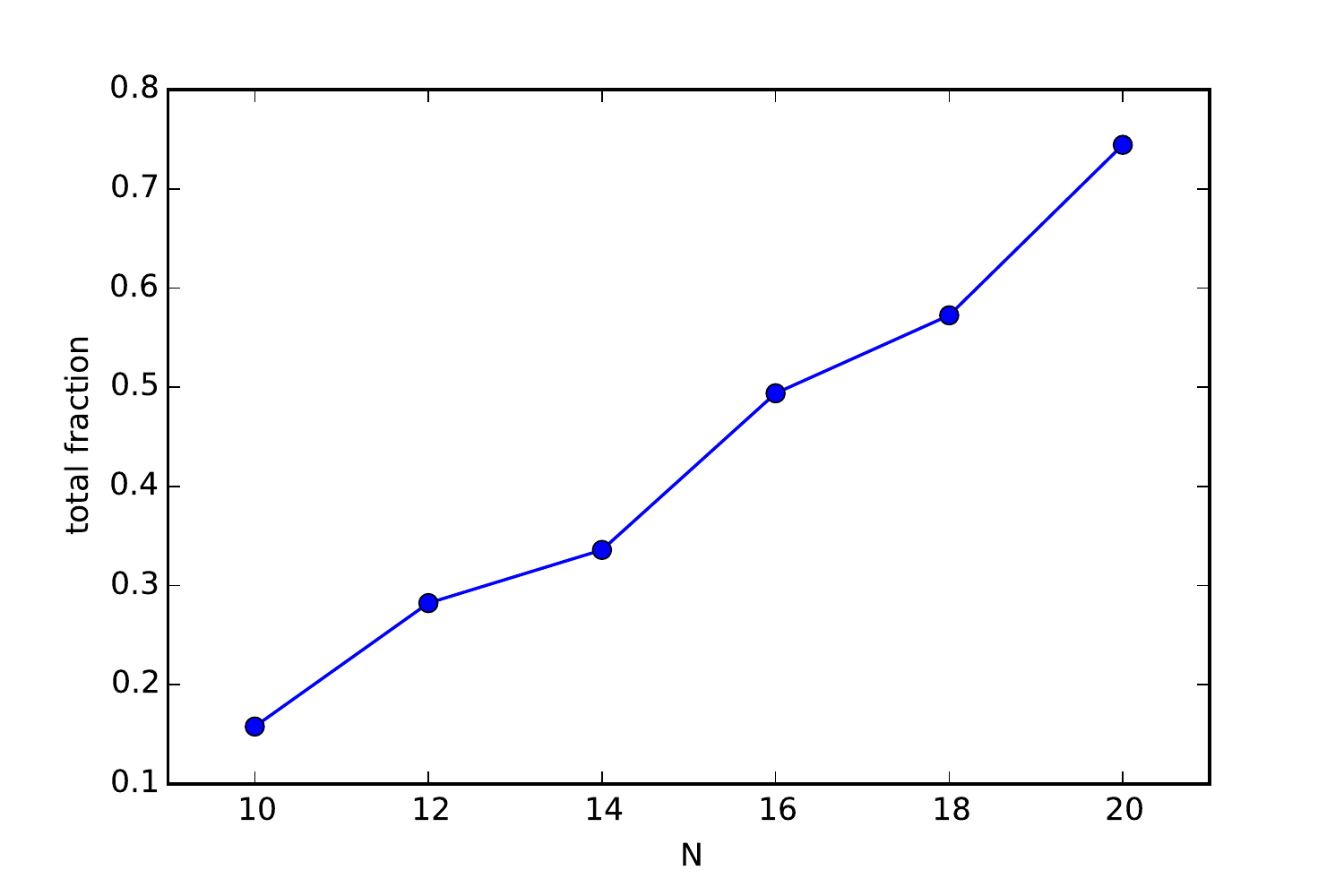}
  \caption{ We plot $\sum_{k}f(N,k)$, where different curves
    correspond to different $N$ and we plot versus $k$.}
  \label{fig:totalpvn}
\end{centering}
\end{figure}

We plot these results in Figures~\ref{fig:pvn},~\ref{fig:pvp},
and~\ref{fig:totalpvn}.  One thing that is clear is that for larger
graphs, the fraction that support a pattern is growing.  

%In fact, when $N=18$, the fraction of cubic graphs that support a
%pattern is more than $57\%$.  We actually get 74\% of the ones with
%$N=20$, but again this is only an approximation since we have only
%done about 2\% of the graphs with $N=20$.  {\bf Upshot:} The data
%seems to strongly suggest that as $N$ grows, the fraction of cubic
%graphs that support a pattern goes to one.

In Figure~\ref{fig:pvn}, each curve corresponds to the fraction of
graphs that support a certain number of patterns.  We see from this
data that each of these functions is increasing as a function of $N$. 

In Figure~\ref{fig:pvp}, each curve corresponds to graphs with a
certain number of vertices.  

Finally, in Figure~\ref{fig:totalpvn} we see that the fraction of
graphs that contain a pattern is monotonically increasing as a
function of the number of vertices. It is nearly linear and based on
the small sample shown in table \ref{tab:patterns}, the plot will
saturate at nearly 100 \% when $N=30$.

%  left column is plotting versus $N$, different curves are different
%    $k$; right column is plotting versus $k$, different curves are
%    different $N$.  Rows are the same data; top row is linear, middle
%    is log, bottom row is plotting $\sum_k p(N,k)$ versus $N$.  Note
%    that the values from $10,\dots,18$ are exact lower bounds, since
%    we did an exhaustive search; the estimate on $N=20$ actually has a
%    sampling error. However, I did a back-of-the-envelope calculation
%    for the 95\% confidence interval and it's about $0.008$, which
%    would barely show up on this graph.  (This is of course assuming
%    that we've chosen a {\bf random} sample of size 9910 from the set
%    of all cubics.) 
%

%\end{document}

\begin{figure}[ht]
\begin{centering}
  \includegraphics[width=.5\textwidth]{{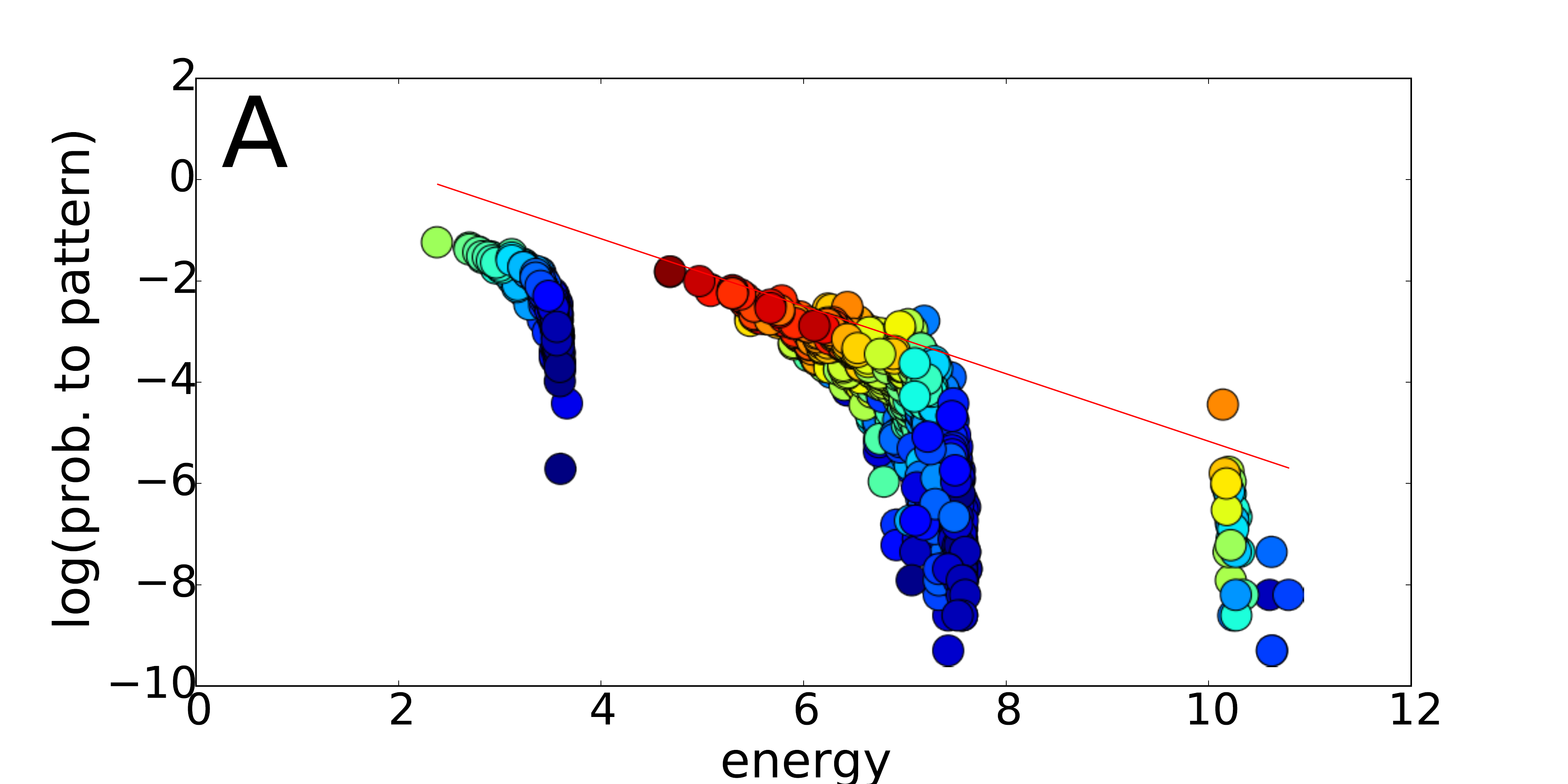}}
  \includegraphics[width=.5\textwidth]{{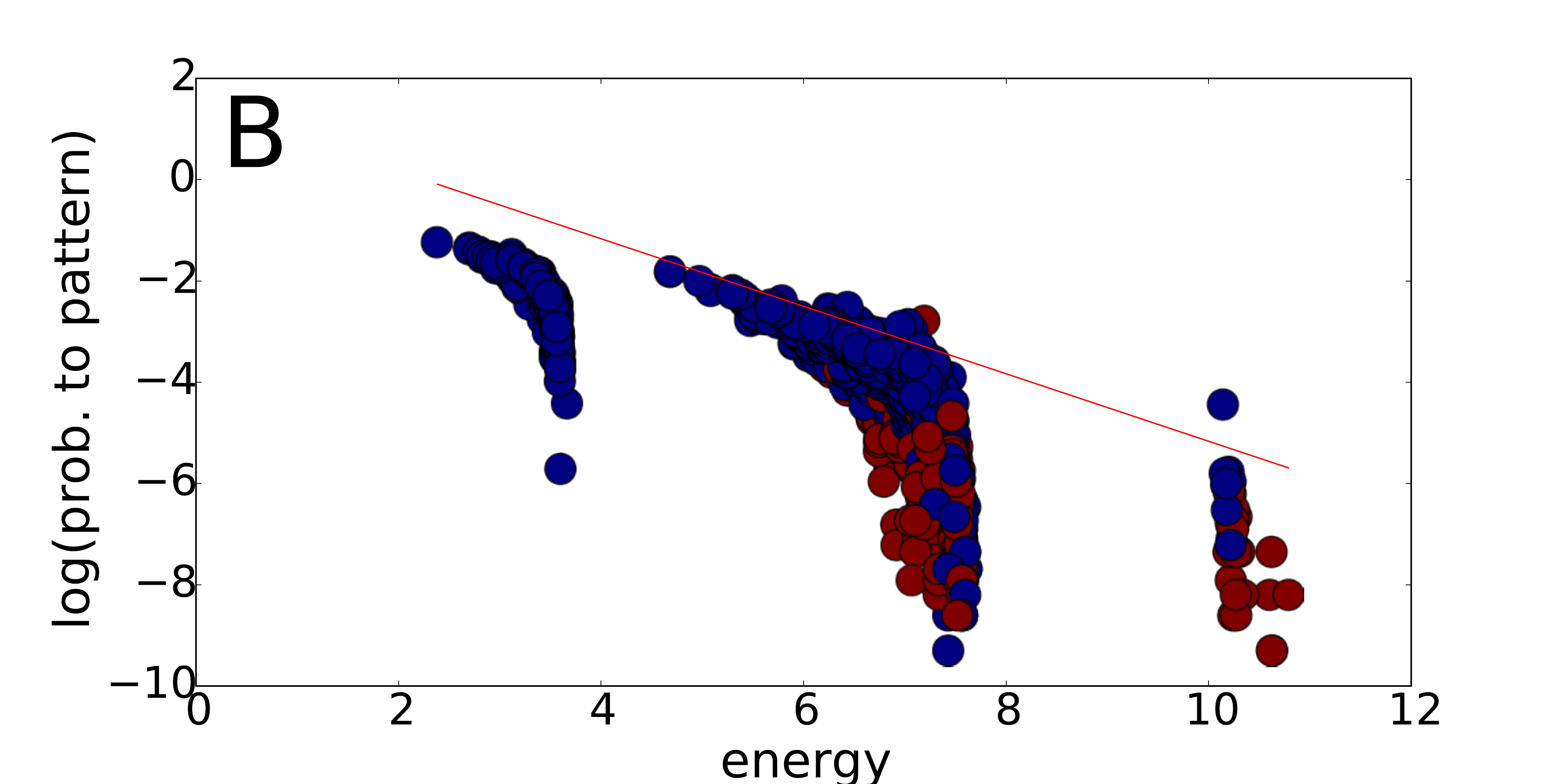}}
  \caption{Plot of all patterns for $N=16$. The same data is plotted
    in both frames with a different coloring scheme.  Each plotted
    circle corresponds to a single pattern; graphs with multiple
    patterns are represented by multiple dots in this plot.  We plot
    the energy of each pattern versus the probability to decay to that
    pattern.  In frame A the dots are colored by the spectral gap,
    which ranges from $0$ (dark blue) to about $0.40$ (red); in frame
    B we are simply coloring due to the presence of a long link: red
    represented a long link solution and blue a short link solution.
    The line corresponds to the function $\exp(-3/2(\mbox{energy}))$
    which seems to be a nice empirical fit (we have no theory for this
    scaling).}
  \label{fig:16}
\end{centering}
\end{figure}

In Figure~\ref{fig:16}, we view all stable patterns for cubic graphs
with $N=16$ vertices.  We are plotting several dimensions of data
here: the energy of the pattern, the probability of it being attained,
the spectral gap of a pattern (the absolute value of the negative
eigenvalue closest to zero), and whether or not the pattern has a long
edge.  The line with slope $-3/2$ has been also added for comparison.

This figure has several salient features: first, we see that the
patterns seem to fall into three clusters; second, inside each of
these clusters, there is a strong (negative) correlation between the
energy of a pattern and the width of its basin of attraction; third,
there is a strong correlation between the width of a basin of a
pattern and its linear stability; fourth, long link solutions are
quite common, but they tend to be higher in energy, harder to hit, and
have smaller spectral gaps.

%\begin{figure}[ht]
%\begin{centering}
%  \includegraphics[width=.5\textwidth]{{N=16longshort.pdf}}
%  \caption{Plot of all patterns for $N=16$, long link versus short link patterns}
%  \label{fig:16ls}
%\end{centering}
%\end{figure}

\begin{figure}[ht]
\begin{centering}
  \includegraphics[width=.5\textwidth]{{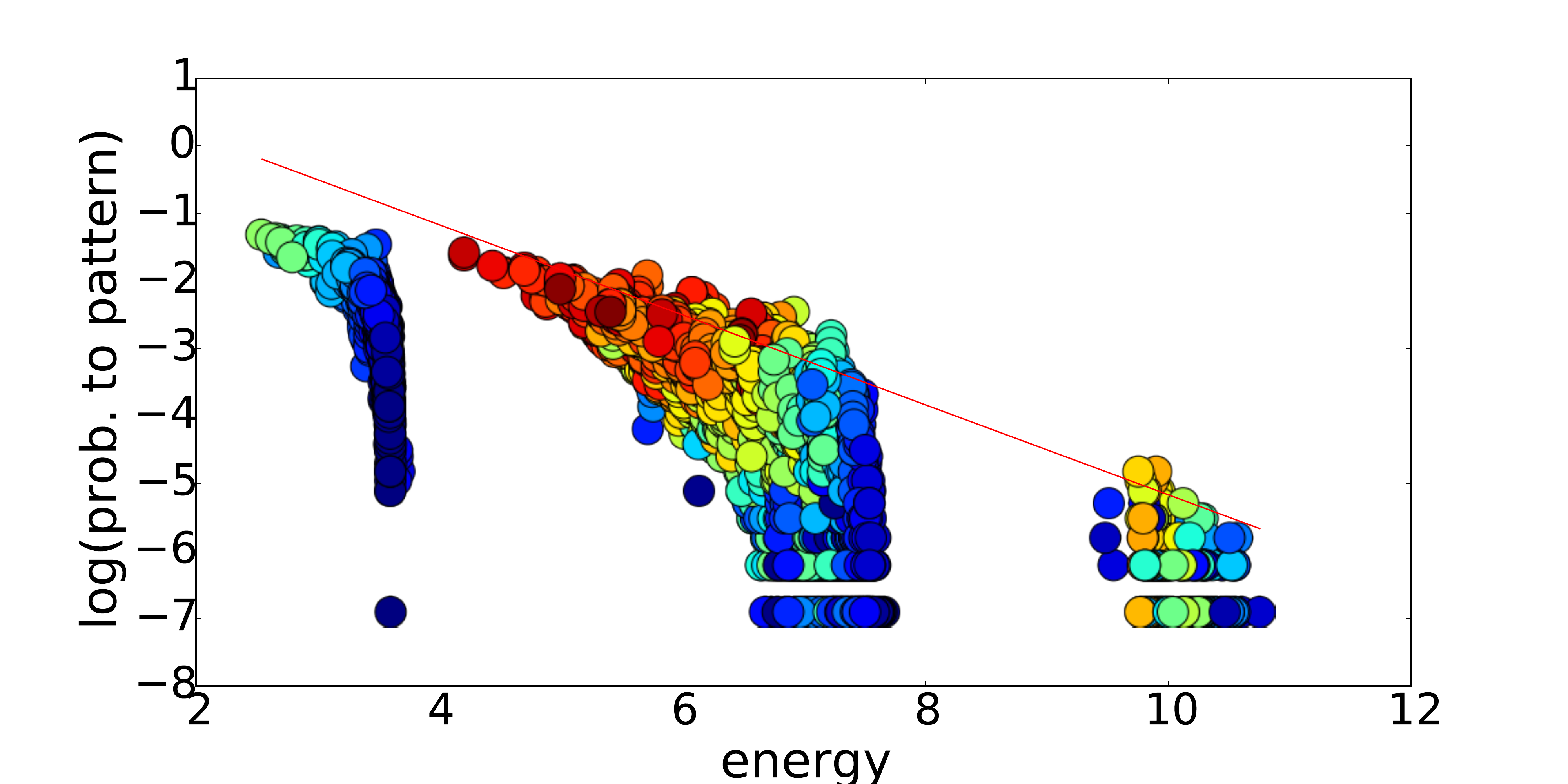}}
  \caption{Plot of all patterns for $N=18$, same scheme used for
    Figure~\ref{fig:16}A.}
  \label{fig:18}
\end{centering}
\end{figure}

\begin{figure}[ht]
\begin{centering}
  \includegraphics[width=.5\textwidth]{{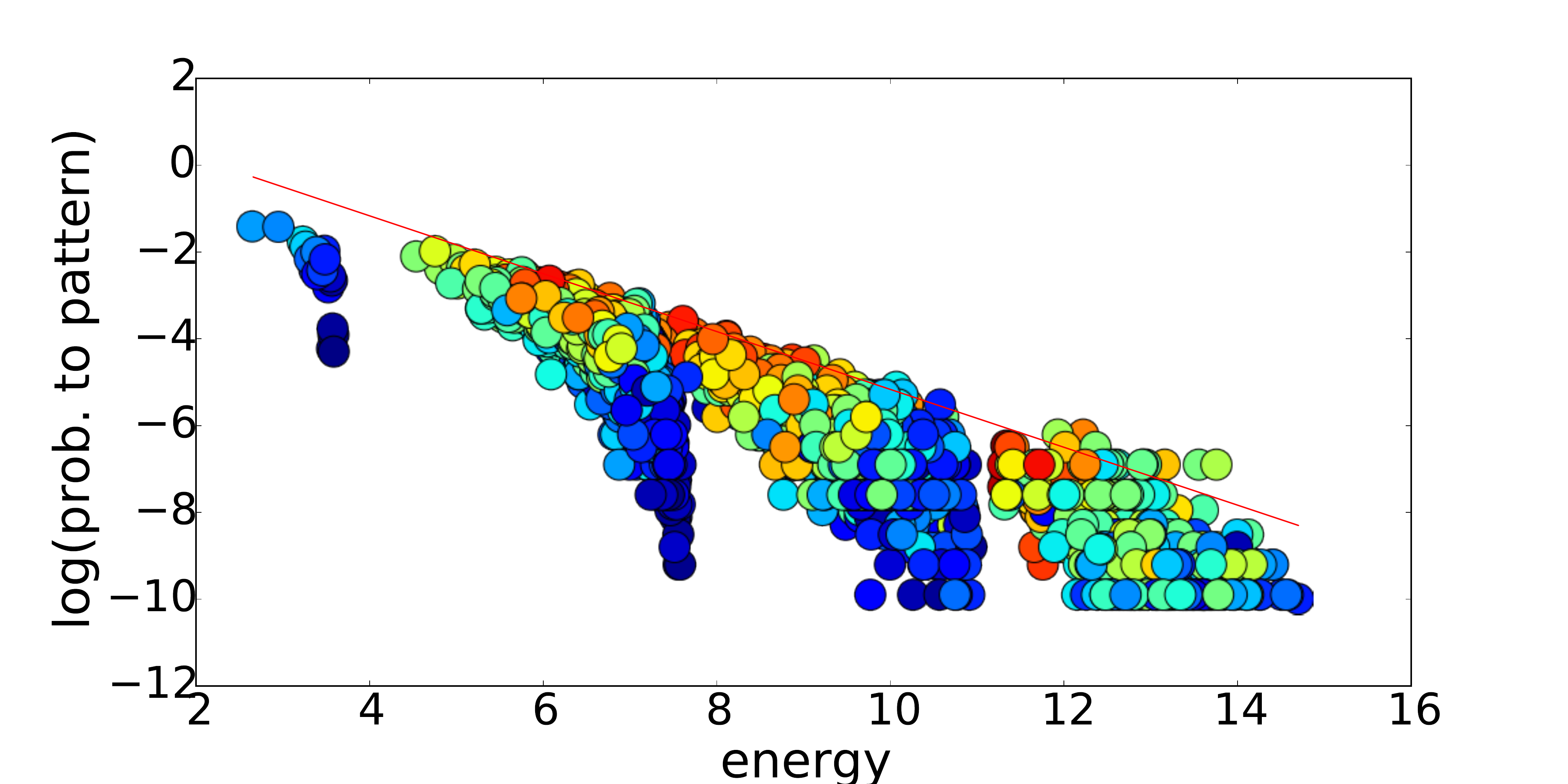}}
  \caption{Plot of all patterns for $N=30$, same scheme used for
    Figure~\ref{fig:16}A.}
  \label{fig:30}
\end{centering}
\end{figure}

In Figures~\ref{fig:18} and~\ref{fig:30}, we plot the same as in
Figure~\ref{fig:16}A, except these are graphs with $N=18, N=30$
vertices respectively.  All of the same patterns apply, with the
exception that there are now four clusters for $N=30$. 

% We also see some ``aliasing'' in the later figures due to finite
% sample size: for example, in the $N=18$ case, we ran each graph for
% $10^7$ samples, and so any time this pattern was found only once,
% this gives a prediction of $10^{-7}$ exactly for the probability to
% hit it.  With larger sample sizes, we would expect to see these
% clusters ``fill out'' as they did in Figure~\ref{fig:16}.

\section{Studying individual graphs}\label{sec:comprehensive}

\subsection{Explaining the clusters}

We observed in several datasets above that the patterns seem to form a
small number of clusters.  For $N=16,18$ there were three, and for $N=30$ there were four.

\begin{figure}[ht]
\begin{centering}
  \includegraphics[width=.4\textwidth]{{{N16.twist1}.pdf}}
  \includegraphics[width=.4\textwidth]{{{N16.twist2}.pdf}}
  \includegraphics[width=.4\textwidth]{{{N16.twist3}.pdf}}
  \caption{$N=16$, one graph from each cluster.}
  \label{fig:16twists}
\end{centering}
\end{figure}

In Figure~\ref{fig:16twists} we plot three graphs, one from each
cluster.  

The manner in which we plot these patterns is as follows: We first
shift each solution so that $\theta_0 = 0$, and then plot the
remaining $\theta$'s at their angle in the plane.  The radius is
chosen to minimize collisions in the picture.  We also color each
vertex by angle, and edges are colored by black if short and red if
long.  Note that each of these graphs has a canonical ``twist'' count,
which is the number of times the pattern wraps around the origin.  It
is very clear in the first two that the winding numbers are one and
two, and in the third one can eventually pick out three distinct
cycles around the origin: one twist is the inner diamond, one is the
trapezoid with two red links, and finally there is a long independent
cycle comprised entirely of short links. Each twist will add a
particular amount of energy to the solution, which we bound below.

For example, assume that we have a twist consisting of $m$ vertices
around the circle; how much energy does this add to the solution?
Assume that we have $m$ vertices located at $0=\theta_0\le \theta_1\le
\dots\le \theta_m=2\pi$.  Writing $\zeta_k = \theta_{k+1}-\theta_k$,
we have that the energy added to the system by this loop is
\begin{equation*}
  F(\zeta) := \sum_{k=0}^{m-1} (1-\cos(\zeta_k)),
\end{equation*}
with the constraints that $\zeta_k\in[0,2\pi], \sum_{k=0}^{m-1}\zeta_k
= 2\pi$.  Using Lagrange multipliers (or symmetry arguments), we see
that the only interior optimizer is at the point $\zeta_k\equiv
2\pi/m$.  The energy of such an equidistant loop is
$m(1-\cos(2\pi/m))$, a function whose maximum value (found at $m=4$)
is approximately $4$.  We see empirically that each loop tends to add
about $3.5$ to the total energy of the pattern. We note that for an
isolated twist to be stable, it has to have at least 5 oscillators
since the phase-difference between connected nodes should be smaller
than $2\pi/4.$ Thus, using $m=5$ for the the energy of the minimal
loop, we get 3.54 which is very close to the energy separating
clusters.  We note that this crude calculation sets an upper bound for
the number of clusters to be roughly $N/5$.  For $N=14$ there are two
clusters with energy difference of about 5.5 (data not shown).

%\todo{Add in some canonical examples from $N=30$}

%\subsection{Comprehensive study for low $N$}

%\subsubsection{Exactly solvable patterns}
%In this subsection, we describe some exactly solvable patterns; that
%is, patterns on graphs that allow us to explicitly write the phase
%relationships. As seen in the datasets above, there is a large number
%of graphs that support patterns, and in this section we are of
%necessity studying only a small subset of these.

\subsection{Double ring graphs and their friends } 

The simplest example of a pattern on a cubic graph is two rings of
length $N/2$ coupled to form a double ring.  To construct this: label
the nodes, $1,2,\ldots, N/2,N/2+1,\ldots, N$. For nodes $1,\ldots,
N/2$, form a ring with nearest neighbor coupling (and 1 and $N/2$
connected) and similarly for nodes $N/2+1,\ldots, 2N.$ Then connect
node $j$ to $N/2+j$ for $j=1,\ldots, N/2.$ The pattern is then
$\theta_j=4\pi(j-1)/N$ for $j=1,\ldots,N/2$ and
$\theta_j=\theta_{j-N/2}$ for $j=N/2,\ldots,N.$ If $N\ge 10$, the
phase-difference between any pairs is less than $\pi/2$ in magnitude,
so that the equilibrium will be stable.  This result proves that there
exists a cubic graph with $N$ nodes that has a pattern for all $N\ge
10.$ The double ring graph has energy $E=N(1-\cos(4\pi/N))$; for
$N=18,E=4.6863$ and this can be seen in figure \ref{fig:16} as the
left most point in the second cluster. (The double ring has two
twists, putting it in the second cluster.)

\begin{figure}[ht]
\begin{center}
\includegraphics[width=.5\textwidth]{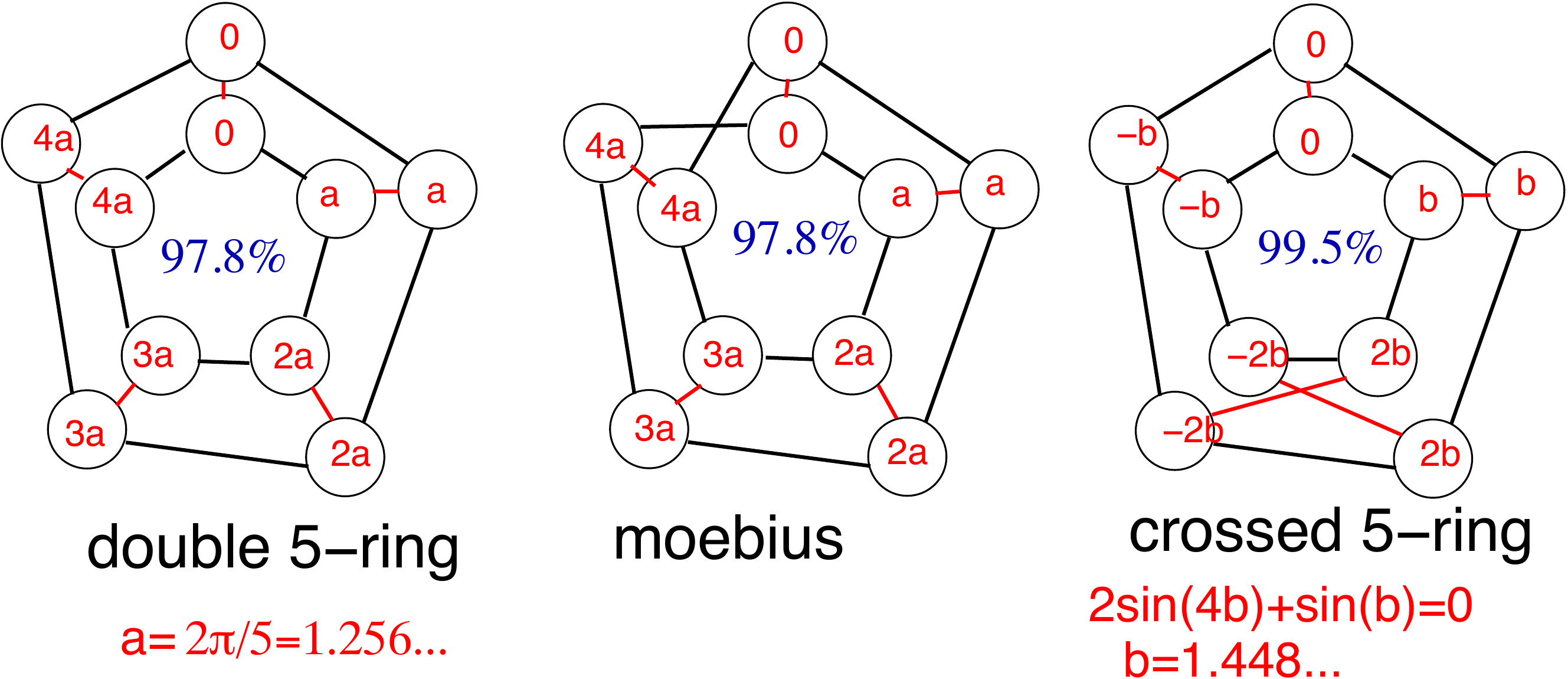}
\end{center}
\caption{ The three cubic graphs of 10 nodes that admit
  patterns. Relative phases are labels $a,b$ and the fraction of
  simulations that lead to synchronization is in blue.  (That is, only
  about 2\% of simulations lead to the pattern shown in the left two
  and 0.5\% in the right-most graph.) The left-most graph is the
  double ring graph for $N=10.$ The middle is a kind of
  M\"{o}bius-strip graph and has the same fixed point and basin width
  as the double ring.  The right-most graph is also like the
  double-ring, but two of the links are crossed.  The value of the
  phase, $b$ is close to, but slightly larger, than $a=2\pi/5$.  The
  phase-difference of $2\pi-4b$ is thus less than $2\pi-4a$ and so is
  still a short link.}
\label{fig:pat10}
\end{figure}

\subsection{$N=10$} 

For ten-node cubic graphs there are only three graphs that have
patterns and each of these has a single pattern. They are illustrated
in figure \ref{fig:pat10}.  The left most is the double-ring graph;
the middle is a M\"{o}bius-strip but has the same phase relations as
the double ring. Both these graphs seem to have a similar basin and
they have exactly the same energy since all phase-differences are the
same.    (Note that the Moebius graph has the
obvious generalization for all $N\ge 10$ similar to the double-ring
construction above.)

The right-most graph is also quite close to the double ring graph, but
the two pairs of nodes have swapped connections.  This graph has a
much smaller basin of attraction (and larger energy) than does the
double ring graph.  Let us consider the analogue of this twisted graph
with more vertices. Let $N=2m.$ Let $b_m$ be the root of the function
\[
g(b)=2\sin((m-1)b)+\sin(b)
\]
that lies between $2\pi/m$ and $2\pi/(m-1).$ It is clear that for
$m>4$ that $g(2\pi/(m-1))>0$ and $g(2\pi/m)<0$ so by continuity such a
root exists. We now define the phases around the twisted circle graph
that consists of two rings of $m$ nodes. We consider $m=2k+1$ and
$m=2k$ separately.  We connect the nodes in each ring exactly as in
the double ring. Next, writing, $m=2k+1$, swap the $k^{th}$ and
$(k+1)^{th}$ connections between the rings to form the twisted
graph. Now assign the phases of the pattern as follows: inner and
outer nodes $j=0,\ldots,k$ are assigned $jb_m$; nodes $j=k+1, \ldots 2
k$ are assigned $-(m-j)b_m.$ This is a phase-locked pattern as long as
$g(b_m)=0.$ (That is the sum of all phase-differences between
connected pairs is zero.)  The possible phase-differences between any
pairs are $0,\pm b_m,\pm (m-1)b_m.$ Our bounds on the root $b_m$
assure that the cosine of each phase-difference is positive as long as
$m>4.$ Thus, that pattern is stable.  For $m$ even, we label the rings
as with the odd case and assign phases on each ring as $0,b_m,\ldots,
(k-1)b_m,-kb_m,-(k-1)b_m,\ldots, -b_m.$ This leads to phase
differences between connected pairs that are either 0, $\pm b_m$, or
$\pm (m-1)b_m.$ Existence and stability are as with the odd case.
This and the previous section show that for any number of nodes
$N\ge10$ there are at least three cubic graphs that have simple stable
patterns that are qualitatively like the traveling wave on a ring of
$N/2$ nodes.

\begin{figure*}[t]
\begin{center}
\includegraphics[width=.9\textwidth]{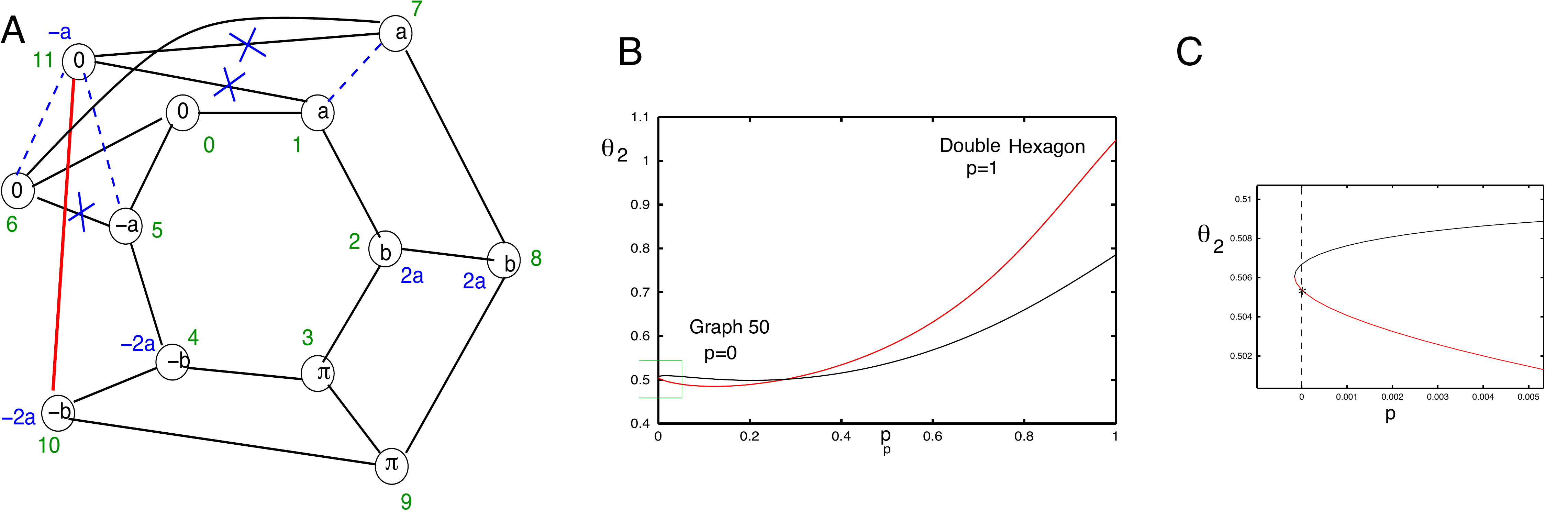}
\end{center}
\caption{A 12-node graph ($G_{50}$) with a long link pattern. (A) The
  graph and its pattern of phases (long link is shown in red) also
  showing that it is similar to the double-ring graph with 12
  nodes. Removal of the X-ed out links and replacing them with the
  dashed links leads to the double ring graph whose phases are shown
  in blue. (B) Numerical continuation of the pattern on the
  double-ring ($p=1$) to the pattern on $G_{50}$ ($p=0$) and an
  expanded view (C).}
\label{fig:g50}
\end{figure*}

\subsection{$N=12$ with long links} \label{sec:N=12long}

For $N=12$ there are 24/88 graphs that support patterns.  $N=12$ is
also the lowest number of nodes for which there are patterns with long
links. Recall long links are those where the phase-difference exceeds
$\pi/2.$ Of the 24 graphs with stable patterns, 4 of them have
patterns with long links and of these 4 graphs, 3 of them have
phase-differences that are very close to $\pi/2$, so the links are
nearly neutral ($\cos(\phi)=0$).  However, one graph has a link where
the phase-difference is quite large and we now discuss this graph (we
will call it $G_{50}$ since was the 50th graph in the database of 12
node cubic graphs referenced in Section~\ref{sec:questions}) as it has
some interesting properties.  It is illustrated in figure
\ref{fig:g50}A and has been drawn in such a way that one can see the
relationship between it and the double 6-ring graph.  First, if the
links drawn with the blue X's are deleted and the dashed blue links
are added instead, we recover the double 6-ring graph and the
corresponding phases are given in the nodes and outside the nodes in
blue, when they differ from those in $G_{50}$.  In this case, the
value of $a$ will be $2\pi/6.$ On the other hand, the phases for
$G_{50}$ are not so symmetric, but are still simple and can be defined
by two numbers, $a$ and $b$.  By looking at the summed sines of the
phase-differences, it can be seen that if there is such a pattern, we
must have $\sin(b-a)=2\sin(a)$ and $\sin(b-a)=\sin(b).$ We can use
Maple or some other symbolic algebra package to see that
\begin{eqnarray*}
a&=&2\asin(1/4)\approx 0.50536,\\
b&=&\pi-\atan\sqrt{15}\approx 1.8235,
\end{eqnarray*}     
so, in particular, $b>\pi/2.$ We see that the long link occurs between
nodes 10 and 11.  Since there are only three connections that are
different between this graph and the double ring, we can start with
the double ring graph pattern and use continuation to homotopy to
$G_{50}$.  Figures \ref{fig:g50}B,C show that the pattern on the double
ring is continuously deformed to the long-link pattern on $G_{50}.$
Note that the limit point or fold occurs very close to $p=0$, showing
that the pattern on $G_{50}$ has an eigenvalue very near zero.  
We also see in the bifurcation diagram of~\ref{fig:g50}C that there is
a saddle-node bifurcation for $p$ a bit less than zero.  This saddle
plays a significant role in the dynamics of Kuramoto oscillators when
noise is added.  For example, the energy difference between the
nearest saddle and the sink is the quantity that governs the mean
escape time away from the sink; more precisely, if we
force~\eqref{eq:K} with white noise of amplitude $\sqrt{2\e}$, then
the mean escape time out of this pattern would be
\begin{equation*}
  \E[\tau] \asymp \exp(\e^{-1} (\Phi_G(\mbox{saddle}) - \Phi_G(\mbox{sink}))).
\end{equation*}
We have not been able to study this systematically for all of the
graphs considered in this paper, but we have seen a strong correlation
in the numerics between three properties for any particular pattern:
(i) the Jacobian at the pattern having a small spectral gap; (ii) the
basin of attraction of the pattern being small; (iii) the potential
height needed to escape the basin of attraction of the pattern.  

In fact, we conjecture that for any pattern, if any of these are
small, then the other two must be.  We can see the correlation between
(i) and (ii) in Figures~\ref{fig:16} and~\ref{fig:30}: there is a
clear pattern that as we move vertically down in the graph, the points
become closer to blue.  Similarly, we see that (ii) and (iii) are
correlated in at least the $G_{50}$ vs. double hex graph, in that the
double hex graph has a much larger basin of attraction.

\begin{figure}[ht]
  \includegraphics[height=4in]{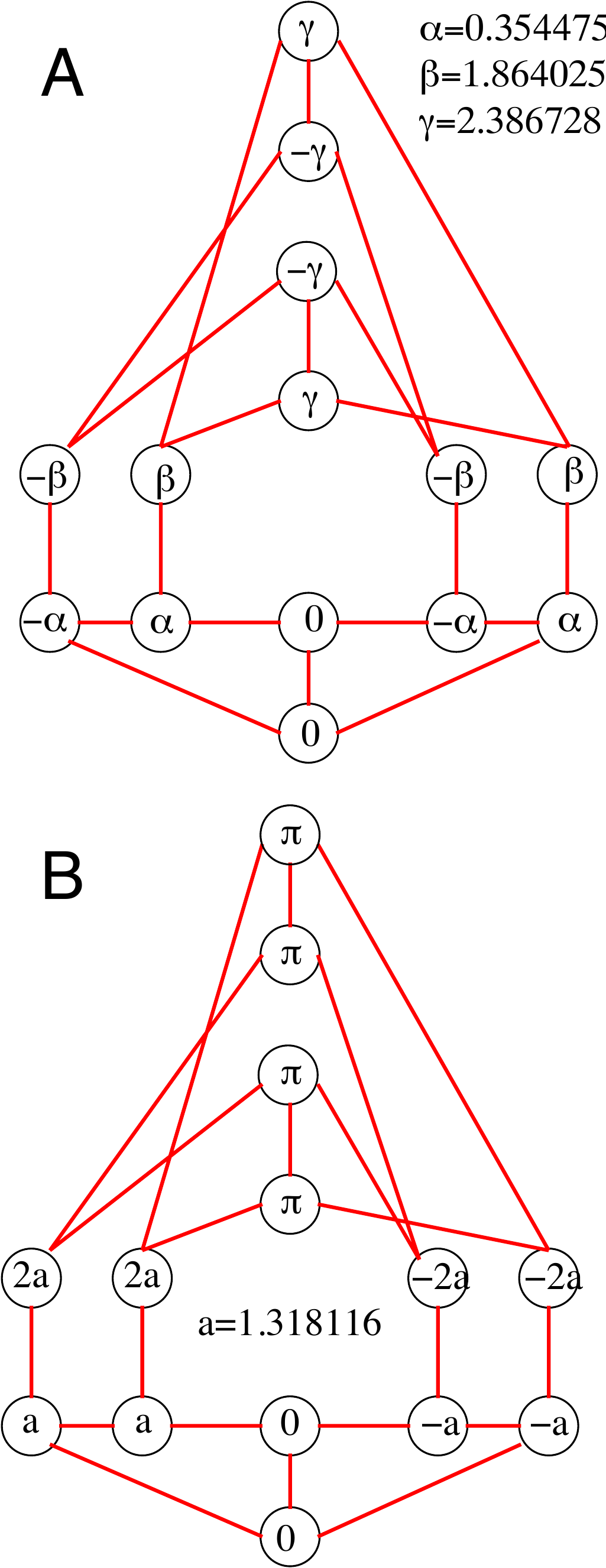}
  \caption{A Graph with two patterns. (A) The higher energy patterns
    depends on the parameters $\alpha,\beta$ with
    $2\gamma+\beta-\alpha=2\pi.$ (Equations for $\alpha,\beta$ are in
    the text.) (B) The low energy pattern depends on the parameter
    $a=\acos(1/4).$ }
\label{fig:twopat}
\end{figure}

\subsection{A graph with multiple patterns.}  $N=12$ is also the
first instance that indicates a graph with more than one possible
pattern; there appear to be two such graphs.  However, the patterns on
these graphs are complex and there is very little difference in their
energies; the exact phases cannot be easily determined analytically,
so we turn to a particular example with $N=14.$ Figure
\ref{fig:twopat} shows such a graph with its two patterns shown in
A,B. The low energy pattern (panel B) depends on only one parameter
that satisfies $2\sin(2a)-\sin(a)=0$ which leads to $a=\acos(1/4).$
The higher energy pattern shown in panel A depends on two parameters,
$\alpha,\beta$ and satisfy:
\begin{eqnarray*}
0&=&\sin(\beta-\alpha)-\sin(2\alpha)-\sin(\alpha)\\
0&=&\sin(\beta-\alpha)-2\sin(3\beta/2-\alpha/2)
\end{eqnarray*}
which can be solved numerically to yield two nontrivial solutions. One
of these solutions leads to an unstable equilibrium, but one is
stable, and this is shown in panel A. The lower energy pattern
($E=7.0$) appears in 1.706\% of the simulations and the higher energy
($E=7.43$) appears in 0.291\% of the simulations; all other
simulations lead to synchrony.

\subsection{High energy patterns.}  Let us define
\begin{eqnarray*}
  E_N &=& 10\left\lfloor\frac{N}{10}\right\rfloor (1-\cos(2\pi/5)), \\
  F_N &=& 10\left\lfloor\frac{N}{10}\right\rfloor(8(1-\cos(\beta^*)) + 4(1-\cos(4\beta^*))),
\end{eqnarray*}
where $\beta^*$ is the root of $2\sin(4x)+\sin(x)$ that lies between
$2\pi/5$ and $\pi/2$.  Then we will show in this section that for
every $N\ge 10$, there is a pattern with energy $E_N$, and for any
$N=10m$ with $m\ge 1$, there is a pattern with energy $F_N$.  Note
that the second statement is always stronger, since
\begin{equation*}
  \frac{E_{10m}}{10m}\approx 6.90983,\quad   \frac{F_{10m}}{10m}\approx 7.49165.
\end{equation*}

The proof for the $F_{10m}$ energy is simpler.  Consider the crossed
5-ring in Figure~\ref{fig:pat10}, where we number the outside vertices
$1,\dots,5$ and the inside vertices $6,\dots,10$, so that
$i\leftrightarrow i+5$.  Consider $m$ disjoint copies of this graph,
with vertices labeled $j.k$, with $j=1,\dots,m$ and $k\in1,\dots,10$.  Now
rewire these graphs as follows: remove the edges $1.1\leftrightarrow
1.6$ and $2.1\leftrightarrow2.6$ and cross them up as
$2.1\leftrightarrow1.6$ and $1.1\leftrightarrow 2.6$.  If $m\ge3$, do
the similar crossing with $2.2\leftrightarrow2.7$ and
$3.2\leftrightarrow3.7$, etc.  Each vertex is still degree three, and
the pattern is still the same, since we have only disconnected and
reconnected vertices with the same angle.  We also see that this
pattern has energy $F_{10m}$.

%The basic idea of the construction is to construct a graph that
%supports $2m$ distinct rings of length 5, with the angles in each ring
%located at multiples of $2\pi/5$.  Each ring will then contribute
%$5(1-\cos(2\pi/5))$ to the total energy of the pattern.  However, we
%still need to show that we can connect these rings in such a manner as
%to obtain a cubic graph.

Now we show that there are stable patterns with energy $E_N$ for all
$N\ge 2$.  We first describe the construction for $N=10m$.  Create
$2m$ 5-ring graphs. These are disjoint and we label their vertices
$A,B,C,D,E.$ Put the 5-wave on each one so that vertex $A$ has phase
0, $B$ has phase $2\pi/5$, etc.  Clearly the energy from each ring is
$5(1-\cos(2\pi/5))$ and the total energy at this point is $E_N$.
Choose some pairing of the $2m$ rings into $m$ pairs, and connect the
$A$ vertices between the two rings in each pair.  Choose a different
pairing of the $2m$ rings, and connect the $B$ vertices, etc.  Go all
the way until the $E$ vertices.  Notice that the connections that came
from the pairings only couple nodes at the same angle, so these new
pairings contribute neither to the vector field~\eqref{eq:K} or the
energy~\eqref{eq:defofPhi}.  Therefore we still have a fixed point with
energy $E_N$.  Also notice that each node has degree three: connected
to two nodes in the original rings, and one more connection in these
pairings.  We can also see that there are choices which give a
connected graph: if we label the original rings $1,\dots, m$ and
$1',\dots, m'$, and choose the $A$ pairing to be $i\leftrightarrow i'$
and the $B$ pairing to be $i \leftrightarrow (i+1)'$, then the entire
graph will be connected at that stage (for example, there is a path in
the rings that goes as $1\to2'\to2\to3'\to3\to\dots$ and the rings are
connected.

One can construct graphs with energy $E_N$ in a similar manner for the
in-between cases, i.e.  $10m < N < 10(m+1)$ where $m\ge 2$.  Let us
write $k=N-10m$, noting that it is even and less than 10.  Let us
write $k/2$ pairs of nodes labeled $A'$, $B'$, using whichever letters
are necessary.  Construct the $10m$ vertex graph from $2m$ rings
described above.  Choose any two pairs of $A$ vertices that were
connected in the last step above, call them $A_1 \leftrightarrow A_2$
and $A_3 \leftrightarrow A_4$.  Break these two connections, and add
the edges
\begin{equation*}
  A'_1\leftrightarrow A_1,A_3,\quad A_2'\leftrightarrow A_2,A_4,\quad A_1'\leftrightarrow A_2'.
\end{equation*}
Note that each of the four $A_i$ vertices have had one edge removed
and one edge reconnected, so are still degree $3$, and note that by
construction $A_i'$ each have degree three.  Thus we now have $10m+2$
vertices in the cubic graph which support the same pattern and energy
$E_m$.  We can repeat this construction as necessary through $B',C'$,
etc.

We see directly that $F_{10}$ is the maximal energy for $N=10$, and we
have not observed any higher energies in the simulations presented in
Section~\ref{sec:global}.  For $N=20$, the highest energy observed by
random search has energy less than 11, yet $E_{20}=13.82$ and
similarly the highest energy we saw for $N=30$ was less than 15, yet
$E_{30}=20.72.$ The reason for this is that our sampling of graphs was
not exhaustive and so there is a low probability that these special
graphs would be chosen.  Furthermore, the actual pattern seems to have
a small basin of attraction.  For $N=20$, we constructed a graph using
the algorithm above that had energy $E_{20}$. We ran simulations with
20000 random initial conditions and the energy $E_{20}$ showed up 3
times. (There were two other patterns that have much lower energies
and occur with a probabilities, 13\% and 0.6\%.)  Thus, these high
energy patterns seem to have a small basin of attraction even when
they do exist. And this actually makes sense: the stable pattern has a
very specific choice of angles and twists around certain cycles, and
it seems unlikely that any initial condition that decayed to it would
have to have the same winding index on all of these cycles.

\begin{figure}[ht]
\begin{centering}
  \includegraphics[width=.45\textwidth]{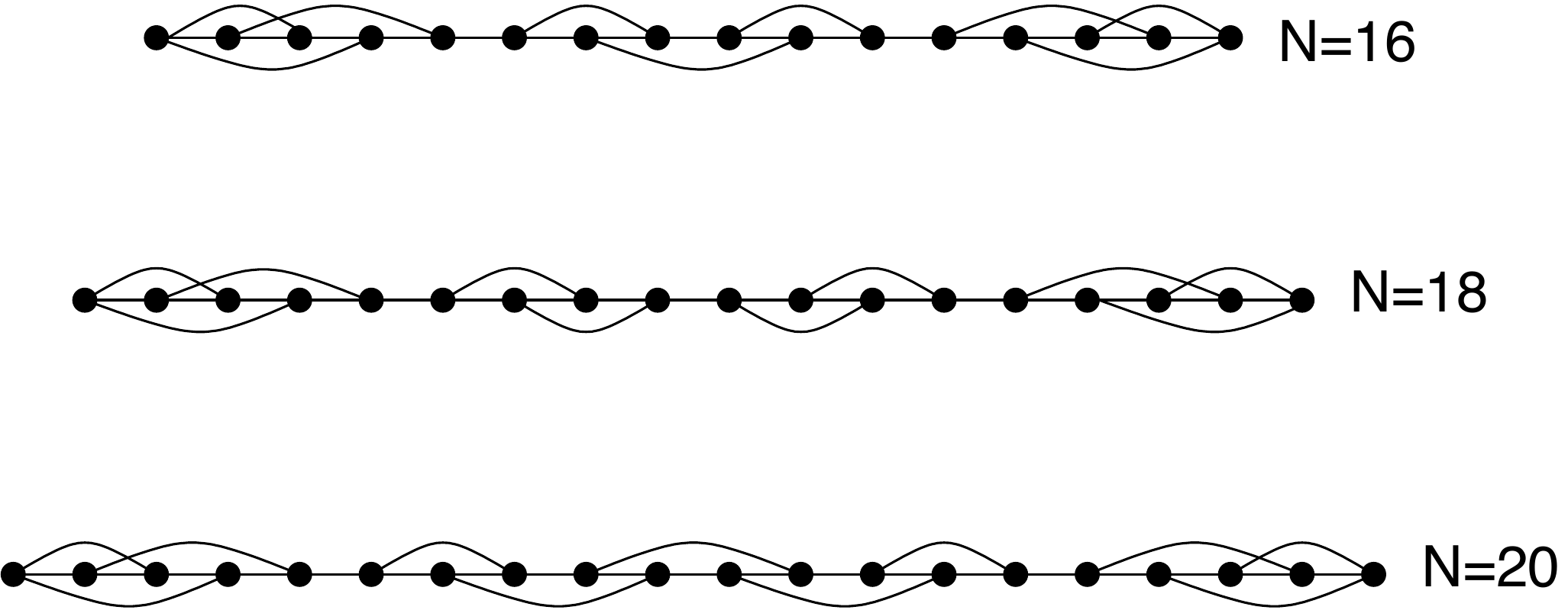}
  \caption{}
  \label{fig:stringy}
\end{centering}
\end{figure}

\subsection{Graphs with no patterns.}
As our statistics have shown, as $N$ gets larger, it seems that there
will be a larger and larger fraction of graphs that have patterns and
this fraction tends to 1. That means, in particular, that graphs with
no patterns are increasingly rare.  We can construct graphs that seem
to have no patterns for any $N$ by starting with the linear chain of
$N$ nodes with nearest neighbor coupling. The two end nodes are
deficient by two edges and the remaining are deficient by 1. The
strategy is to connect the nodes in such a way as to preclude any
cycles that contain 5 or more nodes and such that within those cycles
there are no shortcuts.  That is the shortest path to any node in the
cycle requires going on the cycle.  Figure~\ref{fig:stringy} shows
three graphs for $N=16,18,20.$ The algorithm involves finishing off
the edges first so that there are $N-10$ edges remaining. (Each edge
takes 5 nodes.)  Finish the remaining interior nodes by adding paths
that are either two or three oscillators away.  No cycle contains more
than 4 nodes. (There are longer cycles in the graphs, but there are
shortcuts between nodes in the cycle.) In unpublished work with
Dhaghash Mehta, we have used the methods in~\cite{mehta15} to show
that the three graphs illustrate in the figure have no patterns.  We
conjecture that graphs of this form (which exist for all $N\ge10$)
have no patterns.  If one ignores everything but the nearest neighbor
connections, then it is easy to show that the only phase-differences
are 0,$\pm \pi$ and only 0 is stable.

\section{Conclusions}

We have presented the results of a systematic computation on the
existence and properties of non-trivial phase-locked patterns for the
Kuramoto oscillator on a cubic graph.

The main observations are that cubic graphs can in general support a
wide variety of patterns when there are enough vertices.  We even have
a graph in the dataset that supports twelve distinct patterns (one
example with $N=30$), but graphs that support multiple patterns are
quite common, see Table~\ref{tab:patterns}.  The data also seems to
suggest that as $N$ increases, the probability of a randomly chosen
cubic graph with $N$ vertices not being able to support a pattern
asymptotes to zero.

We have also seen that there are a large number of graphs that support
a pattern with long links, q.v.~Figure~\ref{fig:16}.  In fact, about
one quarter of the patterns observed for $N=16$ are ``long link''
patterns.  In some sense, the search for such patterns is what
motivated the current study, as it is known that the all-to-all
graph~\cite{MS90} and the ring graph~\cite{NON-D} do not have any such
patterns, yet the results of~\cite{Bronski.DeVille.Park.12} show that
it is possible for a long-link pattern to be linearly stable.  In the
event, we have discovered that they are actually quite common.

We have also shown strong numerical evidence that there is a
correlation between the width of a basin of attraction and its linear
stability (in Figures~\ref{fig:16},~\ref{fig:30}, the datapoints get
more blue as we move down in probability).  We have also discovered
that there is a (negative) correlation between this and the energy of
the pattern, but it is somewhat complicated.  In those figures, we see
that all of the patterns form clusters based on their winding number,
and inside each cluster there is a pretty strong negative correlation
between energy and width of the basin.  We also see in some particular
cases a connection between the width of a basin of attraction and its
depth, and conjecture that this relation holds more generally.

However, we see that this study raises more questions than it answers.
It would be interesting to make progress on proving the conjecture
stated in Section~\ref{sec:N=12long}.  Another question would be if we
could determine whether or not a graph would support a pattern based
on some general properties of a graph --- as we mention above, one can
eliminate the possibility for a graph to sustain a pattern when it
lacks certain structures.  But the only analytic technique we know to
show that a graph will support a pattern is the type of hands-on
analysis done in Section~\ref{sec:comprehensive}.

Another interesting question is whether graphs denser than cubic
graphs will have the same propensity for patterns and also for
long-link patterns as we have seen here.  Adding edges can really have
opposing effects when it comes to the stability of a pattern: more
edges means more loops, but it also means more constraints pulling a
pattern in more directions.  One natural direction of future research
is to study more general graphs to try and predict the relationship of
a graphs propensity to support patterns with some gross graph
statistics, e.g. average degree, average clustering coefficient, etc.

\section{Acknowledgments}
We would like to thank Dhagash Mehta and Jonathan Hauenstein for
computing the all fixed points for the graphs in figure
\ref{fig:stringy}. We would also like to thank summer student Ethan
Levien for some early computations of non-synchronized attractors in
the small cubic graphs. B.E. was supported by NSF grant DMS 1219753;
L. D. was supported by the National Aeronautics and Space
Administration (NASA) through the NASA Astrobiology Institute under
Cooperative Agreement Number NNA13AA91A issued through the Science
Mission Directorate.

%\bibliographystyle{amsalpha} 
%\bibliography{kuram}

\providecommand{\bysame}{\leavevmode\hbox to3em{\hrulefill}\thinspace}
\providecommand{\MR}{\relax\ifhmode\unskip\space\fi MR }
\providecommand{\MRhref}[2]{%
  \href{http://www.ams.org/mathscinet-getitem?mr=#1}{#2}
}
\providecommand{\href}[2]{#2}

\end{document}